\documentclass[smallextended]{svjour3}       % onecolumn (second format)

\pdfoutput=1
\usepackage{algorithm}
\usepackage{algpseudocode}
\usepackage{algcompatible}

\usepackage{color}
\usepackage[justification=centering]{caption}

\usepackage{graphicx}
\usepackage{amssymb,amsmath,bm,mathrsfs}
\newcommand{\diag}{\mbox{\rm diag}}

\usepackage{mathtools}

\DeclarePairedDelimiter\floor{\lfloor}{\rfloor}

\def\a{\alpha}
\def\b{\beta}
\def\de{\delta}
\def\De{\Delta}
\def\ga{\gamma}
\def\Si{\Sigma}
\def\si{\sigma}

\def\ze{\zeta}
\def\om{\omega}
\def\leq{\leqslant}
\def\rar{\rightarrow}

\def\R{\hro{R}}
\def\C{\hro{C}}
\def\hro{\mathbb}
\def\N{\hro{N}}

\def\tA{\tilde{A}}

\def\tB{\tilde{B}}
\def\hB{\hat{B}}
\def\tN{\tilde{N}}

\def\tf{\tilde{f}}

\def\lb{\lambda}
\def\cP{{\cal P}}
\def\cQ{{\cal Q}}
\def\fr{\frac}
\def\frQ{\mathfrak{Q}}

\def\cU{{\cal U}}
\def\cV{{\cal V}}

\def\tcU{{\tilde{\cal U}}}
\def\cZ{{\cal Z}}

\def\ddt{\fr{\mathrm{d}}{\mathrm{d}t}}

\def\tU{\tilde{U}}

\def\tcQ{{\tilde{\cal Q}}}
\def\tcP{{\tilde{\cal P}}}
\def\hcP{\widehat{\cP}}
\def\hcQ{\widehat{\cQ}}

\def\vphi{\varphi}

\newcommand{\n}[1]{\left\lVert#1\right\rVert}

\newcommand{\m}[1]{
\begin{bmatrix}
 #1
\end{bmatrix}
}

\renewcommand{\pm}[1]{
\begin{matrix}
 #1
\end{matrix}
}

\def\ES{
	\begin{bmatrix}
		J^E  & 0        & 0       & 0 \\
		0       & N^E_{2} & 0       & 0 \\
		0       & 0        & N^E_{3}& 0 \\
		0       & 0        & 0       & N^E_{4} 
	\end{bmatrix}
}

\def\AS{
	\begin{bmatrix}
		A_{1}  & 0       & 0          &  0        \\
		0           & J^A  & 0          & 0      \\
		0           &   0                & N^A_{3}      & 0      \\
		0           &   0                & 0                    & N^A_{4}      
	\end{bmatrix}
}
\def\BS{
	\begin{bmatrix}
		B_{1}  & 0       & 0        & 0 \\
		0          & B_{2}       & 0        & 0 \\
		0          & 0                  & J^B  & 0  \\
		0          & 0                 &   0                & N^B_{4}  \\ 
	\end{bmatrix}
}

\def\be{\begin{equation}}
\def\ee{\end{equation}}         
\newcommand{\ben}{\begin{eqnarray}}
\newcommand{\een}{\end{eqnarray}}
\newcommand{\bens}{\begin{eqnarray*}}
\newcommand{\eens}{\end{eqnarray*}}
\def\bc{\begin{cases}}
\def\ec{\end{cases}}
\newcommand{\bsq}{\begin{subequations}}
\newcommand{\esq}{\end{subequations}}

\begin{document}

%%%%%%%%%%%%%%%%%
\title{Spectral Characterizations of Solvability and Stability \\ for Delay Differential-Algebraic Equations}
%\subtitle{Spectral Characterizations of Solvability and Stability for DDAEs}
\titlerunning{Spectral Characterizations of Solvability and Stability for DDAEs}% Part of RIGHT running header
\author{ {\sc Phi Ha \thanks{This research is funded by the Nafosted foundation of Vietnam under the project number 101.01-2017.302}
}}
%\authorrunning{Short author list}% Part of LEFT running header
\institute{Phi Ha \at Hanoi University of Science, VNU \\ 	Nguyen Trai Street 334, Thanh Xuan, Hanoi, Vietnam\\ 
	       \email{haphi.hus@vnu.edu.vn}
          }
\date{Received: date / Accepted: date}      

\maketitle

%%%%%%%%%%%%%%%%%abstract style
%Two grouping braces are necessary in abstract environment
%first argument contains abstract text; second argument contains keywords
%text
\begin{abstract}
{The solvability and stability analysis of linear time invariant systems of delay differential-algebraic equations (DDAEs) is analyzed. The behavior approach is applied to DDAEs in order to establish characterizations of their solvability in terms of spectral conditions. Furthermore, examples are delivered to demonstrate that the eigenvalue-based approach to analyze the exponential stability of dynamical systems is only valid for a special class of DDAEs, namely \emph{non-advanced}. Then, a new concept of weak stability is proposed and studied for DDAEs whose matrix coefficients pairwise commute.
}
\end{abstract}


%%%%%%%%%%%%%%section 1%%%%%%%%%
\section{Introduction and Preliminaries}\label{introduction}

%%%%%%%%%%%%%%%%%%%%%%%%%%%%%

Our focus in the present paper is on the solvalbility analysis of linear, time invariant {\it delay dif\-fe\-ren\-tial-al\-ge\-braic equations (DDAEs)} of the form
\begin{equation}\label{eq1.1}
  E\dot{x}(t) = A x(t) + Bx(t-\tau) + f(t), \ \mbox{ for all } t\in [0,\infty),
\end{equation}
and the stability analysis of its associated homogeneous system 
\be\label{homDDAE}
E\dot{x}(t) = A x(t) + Bx(t-\tau), \ \mbox{ for all } t\in [0,\infty),
\ee
where $E$, $A$, $B\in \R^{\ell,n}$, $x:[-\tau,\infty) \rar \R^n$, $f:[0,\infty)\to \R^{\ell}$, and $\tau>0$ is a constant delay. 

%\newpage 

DDAEs of the form \eqref{eq1.1} can be considered as a general combination of two important classes of dynamical systems, namely \emph{differential-algebraic equations} (DAEs)
\begin{equation}\label{eq1.2}
 E \dot{x}(t) = A x(t) + f(t),
\end{equation}
where the matrix $E$ is allowed to be either non-square ($\ell\not= n$) or singular ($\det E=0$), and \emph{delay-differential equations} (DDEs)
\begin{equation}\label{eq1.3}
 \dot{x}(t) = A x(t) + B x(t-\tau) + f(t). 
\end{equation}
Due to the broad range of applications of both DAEs and DDEs, DDAEs have been arisen in various applications, see \cite{AscP95,Cam80,HalL93,ShaG06,ZhuP97} and the references there in. From the theoretical viewpoint, since DDAEs are matrix delay differential equations coupled with matrix difference equations, the study for such systems is much more complicated than that for standard DDEs or DAEs. The dynamics of DDAEs, therefore, has been strongly enriched, and many interesting properties, which occur neither for DAEs nor for DDEs, have been observed \cite{Cam95c,DuLMT13,HaM12,HaM16}. Due to these reasons, recently more and more attention has been devoted to DDAEs, \cite{CamL09,Fri02,HaM12,HaM16,Mic11,ShaG06,TiaYK11,LinT15}.

For both DAEs and DDEs, the structure of their matrix coefficients always play a very important role in the solvability/stability analysis. 
It is well known, that for DDEs of the form \eqref{eq1.3}, stability properties of the solution are closely related to spectral conditions of the matrix triple $(I,A,B)$, see \cite{HalL93}. From the DAE side, not only the stability of \eqref{eq1.2} depends on spectral conditions of the matrix pencil $\lb E-A$ 
but also the solvability is connected to the regularity of this pencil, see Definition \ref{regularity} below. Consequently, both the solvability and the stability of DDAEs are usually discussed under the regularity assumption of this pencil, see e.g., \cite{AscP95,Cam80,CamL09,Fri02,Mic11,ShaG06,TiaYK11,ZhuP97} and the references there in. 
Only a few results for DDAEs have been achieved when the pencil $\lb E-A$ is singular, see e.g., \cite{Cam95c,HaM12,HaM16}. 
This paper aims to fill in this gap.

The short outline of this work is as follows. After some notations and auxiliary lemmas, in Section \ref{sec 2}, 
we construct the condensed form for couples of matrix polynomials. This form allows us to study the solvability of system \eqref{eq1.1} via the structure of the matrix triple $(E,A,B)$. The existence and uniqueness of the solution to \eqref{eq1.1} is, therefore, linked with spectral conditions of \eqref{eq1.1}. Moreover, the presented approach allows further theoretical investigation on more general systems, for example, systems with high order derivatives of both $x(t)$ and $x(t-\tau)$, or systems with either no solution or multiple solutions, etc. 
In Section 3 we first demonstrate that the eigenvalue-based approach is not always suitable to analyze the stability of high index DDAEs, and therefore, either an index reduction procedure or a new stability concepts must be considered. In this section, a new concept of weakly exponential stability is proposed. Finally, in Section 4 we analyze both the solvability and stability of DDAEs whose matrix coefficients pairwise commute. 

In the following we denote by $\N$ ($\N_0$) the set of natural numbers (including $0$), by
$\R$ ($\C$) the set of real (complex) numbers and $\C_{-}:=\{\lb \in \C \ | \ Re \lb <0\}$. 
By $\|\cdot\|$ we denote a norm in $\R^n$, by $\R^{\ell,n}$ 
the real matrices of size $\ell \times n$ and by $I$ ($I_n$) the identity matrix (of size $n \times n$).
As usual $x^{(j)}$ is the $j$-th derivative of a function $x$.

For $ 0 \leq p \leq \infty$, the set $C^p([-\tau,0],\R^n)$ denotes the space of $p$-times continuously differentiable functions from 
$[-\tau,0]$ to $\R^n$. These spaces are equipped with the norm $\|\cdot\|_{C^p}$ defined by 
$\|\vphi\|_{C^p} := \sum_{i=0}^{p} \underset{t\in[-\tau,0]}{\sup} \|\vphi^{(i)}(t)\|$
to form a Banach space. For $p=0$, we adopt the notation $C([-\tau,0],\R^n)$ with the norm 
$\|\cdot\|_{\infty}:=\|\cdot\|_{C^0}$. 
%Furthermore, let
%%
%\[
% C^{\infty}_b([-\tau,0],\R^n):= \{\vphi \in C^{\infty}([-\tau,0],\R^n)| \ \sup_{0 \leq i \leq \infty}\sup_{t\in[-\tau,0]} \|\vphi^{(i)}(t)\| < \infty \},
%\]
%%
%be equipped with the norm
%%
%\[
% \|\vphi\|_{C^{\infty}} := \sup_{0 \leq i \leq \infty} \sup_{t\in[-\tau,0]} \|\vphi^{(i)}(t)\|,
%\]
%%
%to form the Banach space $(C^{\infty}_b([-\tau,0],\R^n),\|\cdot\|_{C^{\infty}})$ \cite{}.

To achieve uniqueness of solutions, analogous to the theory of DDEs, for DDAEs of the form \eqref{eq1.1} one typically has to prescribe an initial function, which takes the form
\begin{equation}\label{eq1.4}
x|_{[-\tau,0]}=\vphi: [-\tau,0] \rightarrow \R^{n}.
\end{equation}
For the DAE \eqref{eq1.2} (resp. the DDE \eqref{eq1.3}), one frequently uses the concept of classical solutions, i. e., functions
which are continuously differentiable and satisfy \eqref{eq1.2} (resp. \eqref{eq1.3}) pointwise, see e.g. \cite{BelZ03,BreCP96}. 
However, there is no clear reason why $E(0)\dot{x}(0)$ which arises in \eqref{eq1.1} should be equal to $E(0)\dot{\vphi}(0^-)$.
Furthermore, for DDAEs, it has been observed in \cite{BakPT02,Cam80,GugH07} that a discontinuity of $\dot{x}$ at $t=0$
may propagate with time, and typically $\dot{x}$ is discontinuous at every point $j\tau, \ j\in\N_0$.
To deal with this property of DDAEs, we use the following solution concept. % for \eqref{eq1.1}.
\begin{definition}\label{solution}
A function $x:[-\tau,\infty)\rar\R^n$ is called a \emph{piecewise differentiable solution} of \eqref{eq1.1}, if $Ex$ is piecewise continuously 
differentiable, $x$ is continuous and satisfies \eqref{eq1.1} at every  $t\in [0,\infty) \setminus \underset{j\in \N_0}{\cup} \{j\tau\}$.
\end{definition}

Throughout this paper whenever we speak of a solution, we mean a piecewise differentiable solution. Notice that, like DAEs, DDAEs are not solvable for arbitrary initial conditions, but they have to obey certain consistency conditions.
\begin{definition}\label{consistency} An initial function $\vphi$ is called \emph{consistent} with \eqref{eq1.1} if the associated initial value problem (IVP) \eqref{eq1.1}, \eqref{eq1.4} has at least one solution.
System \eqref{eq1.1} is called \emph{solvable} (resp. \emph{regular}) if for every consistent initial function $\vphi$,
the associated IVP \eqref{eq1.1}, \eqref{eq1.4} has a solution (resp. has a unique solution).
\end{definition}

\begin{definition}\label{regularity} Consider the DDAE \eqref{eq1.1}. The matrix triple $(E,A,B)$ is called \emph{regular} if the (two variable) \emph{characteristic polynomial} $\mathfrak{P}(\lb,\om):=\det(\lb E - A - \om B)$ is not identically zero. 
If, in addition, $B=0$ we say that the matrix pair $(E,A)$ (or the pencil $\lb E-A$) is regular.  
The sets $\si(E,A,B):= \{\lb \in \C \ | \ \det(\lb E -A-e^{-\lb \tau}B) = 0\}$, $\rho(E,A,B)=\C\setminus \si(E,A,B)$ 
are called the \emph{spectrum} and the \emph{resolvent set} of \eqref{eq1.1}, respectively. 
\end{definition}

%%%%%%%%%%%%%section 2%%%%%%%%%
%\newpage 
\section{Solvability analysis of linear DDAEs}\label{sec 2}

%%%%%%%%%%%%%%%%%%%%%%%%%%%%%
There is one well-known fact in the theory of DAEs, see e.g. \cite{BreCP96,KunM06}, that the DAE \eqref{eq1.2} is uniquely solvable if and only if the matrix pair $(E,A)$ is regular. If this is the case, one can write down explicitly the solution's formula for \eqref{eq1.1}.
Because of this reason, previous works on DDAEs usually consider the case where the pair $(E,A)$ is regular. 
However, in general system \eqref{eq1.1} can be uniquely solvable without the regularity of the pair $(E,A)$, for example the system 
\[
\m{1 & 0 \\ 0 & 0} \m{\dot{x}_1(t) \\ \dot{x}_2(t)} = \m{0 & 0 \\ 0 & 1} \m{{x}_1(t-\tau) \\ {x}_2(t-\tau)} + \m{f_1(t) \\ f_2(t)}, \ \mbox{ for all } t\in [0,\infty).
\]
Therefore, one may expect, that the necessary and sufficient condition for the existence and uniqueness of a solution to \eqref{eq1.1} is the regularity of the matrix triple $(E,A,B)$. This section aims to provide an answer for this question, which has only been rarely considered in literature, \cite{Cam95c}. Rather than working directly with three matrices, we perform the behavior approach \cite{PolW98}, in order to rewrite the considered system in a new form that contains only two matrix polynomials. Then, we analyze the structure of these two matrix polynomials in order to read off the solvability of \eqref{eq1.1}. Notice that, even though the behavior approach has been extensively studied for both DAEs and DDEs \cite{Glu02,PolW98}, it has not been considered for DDAEs. Furthermore, in order to perform this method, it is needed to consider a solution $x$ within the space of infinitely smooth functions, i.e. $C^{\infty}([0,\infty),\R^{n})$. Thus, only within this section, we assume that $x|_{[0,\infty)} \in C^{\infty}([0,\infty),\R^{n})$.

For the sake of completeness, we study the general system
\begin{equation}\label{eq2,2}
  A_k x^{(k)}(t) \!+\! \dots \!+\! A_0 x(t) \!=\! B_{\kappa} x^{(\kappa)}(t-\tau) \!+\! \dots \!+\! B_{0} x(t-\tau) \!+\! f(t), \mbox{ for all } t\in [0,\infty),
\end{equation}
where the coefficients satisfy $A_i \in \R^{\ell,n}$, $B_j \in \R^{\ell,n}$, $i=k,\dots,0$, 
$j=\kappa,\dots,0$.
Let $\cP(\lb):= A_k \lb^k + \dots + A_1 \lb + A_0$, and $\cQ(\lb):= B_{\kappa} \lb^{\kappa} + \dots + B_1 \lb + B_{0}$,
we first rewrite system \eqref{eq2,2} in the behavior form
\begin{equation}\label{eq2.3}
 \cP(\fr{d}{dt}) x(t)=\cQ(\fr{d}{dt})x(t-\tau) + f(t), \ \mbox{ for all } t\in [0,\infty).
\end{equation}
We notice that the characteristic polynomial is $\mathfrak{P}(\lb,\om):=\det(\cP(\lb) - \om \cQ(\lb))$.
\noindent A matrix polynomial ${\cal M} \in \R[\xi]^{n\times n}$ is called \emph{unimodular} if $\det({\cal M})$ is a non-zero constant.
Similar to the singular value decomposition of a matrix, there exists a Smith canonical form for matrix polynomials. 
\begin{proposition}(Smith canonical form \cite{PolW98}) Let $\cP \in \R[\xi]^{\ell\times n}$. Then there exist unimodular matrix polynomials 
$\cU \in \R[\xi]^{\ell\times \ell}$ and $\cV \in \R[\xi]^{n \times n}$ such that
\[ \cU \cP \cV = \m{\diag(p_1,p_2,\dots,p_r) & 0^{r \times (n-r)} \\ 
0^{(\ell-r)\times r} & 0^{(\ell-r)\times (n-r)}} \]
with $p_1,p_2,\dots,p_r \in \R^{1\times 1}[\xi]$ are monic, i.e., their leading coefficients 
are equal to 1, and $p_k$ divides $p_{k+1}$ for $k =1,2,\dots,r$. 
\end{proposition}

Now we construct the condensed form of matrix polynomial pairs in the next theorem.
\begin{theorem}\label{thm2.1} For any matrix polynomial pair $(\cP,\cQ) \in 
(\R[\xi]^{\ell\times n})^2$, there exist unimodular matrix polynomials $\cU \in \R[\xi]^{\ell\times \ell}$, 
$\cV \in \R[\xi]^{n \times n}$ such that
\begin{equation}\label{staircase-form}
 \cU \cP \cV \!=\! 
 \left[ \begin{array}{cc|cccc}
 \Si_{P}  & 0 & *  & *   & \dots  & * \\ \hline
          &   & 0  & *   & \dots  & *  \\
          &   &    & 0   & \dots  & *  \\
          &   &    &     & \ddots & \vdots  \\  
          &   &    &     &        & 0        \\  \hline 
          &   & 0  & *   & \dots  & *  \\ 
          &   &    & 0   & \dots  & *  \\               
          &   &    &     & \ddots & \vdots  \\  
          &   &    &     &        & 0       
 \end{array} \right], \ \cU \cQ \cV \!=\! 
  \left[ \begin{array}{cc|cccc}
 \cQ_{11} &  \cQ_{12} & *  & *   & \dots  & * \\ \hline
          &   & \Si_{q}  & *   & \dots  & *  \\
          &   &    &\Si_{q-1}& \dots  & *  \\
          &   &    &     & \ddots & \vdots  \\  
          &   &    &     &        & \Si_{1}        \\  \hline 
          &   & 0  & *   & \dots  & *  \\ 
          &   &    & 0   & \dots  & *  \\               
          &   &    &     & \ddots & \vdots  \\  
          &   &    &     &        & 0 
\end{array} \right],
\end{equation}
where 
\begin{equation}\label{eq2.2}
 \Si_{j} = 
 \m{ p_{j,1} &         & \\
             & \ddots  & \\
             &         & p_{j,r_j} 
 }
\end{equation}
with monic polynomials $p_{j,1}$, $p_{j,2},\dots,p_{j,r_j}$ on the main diagonal. 
\end{theorem}
\proof
The proof is obtained by a constructive way based on the condensed form approach proposed in \cite{ByeGM97,ByeKM97}. Consider the recursive procedure:

\noindent \textbf{Initial:} Let $(\cP_1,\cQ_1)=(\cP,\cQ) \in (\R[\xi]^{\ell\times n})^2$ 
and set $i=1$. 

\noindent \textbf{Step 1.} Letting $\cU_1$, $\cV_1$ be unimodular matrix polynomials that produce the Smith form of 
$\cP_1$ and partitioning $\cU_1 \cQ_1 \cV_1$ conformably, we get
\[
\cP_2 :=  \cU_1 \cP_1 \cV_1 = \m{\widetilde{\Si} & 0 \\ 0 & 0}, \ 
\cQ_2 := \cU_1 \cQ_1 \cV_1 = \m{\cQ_{11} & \cQ_{12} \\ \cQ_{21} & \cQ_{22}}. 
\]
\textbf{If} $\m{\cQ_{21}&\cQ_{22}} = 0$ then \textbf{Stops}, \textbf{otherwise} proceed to Step 2.

\noindent \textbf{Step 2.} Now we consider the Smith canonical form of $\m{\cQ_{21}&\cQ_{22}}$. Let $\tcU_2$, $\cV_2$ be unimodular matrix polynomials such that
\[
  \tcU_2 \m{\cQ_{21} & \cQ_{22}} \cV_2 = \m{0 & \Si_i \\ 0 & 0},
  \qquad \pm{a_i \mbox{ rows }\\ v_i \mbox{ rows }}
\]
where $\Si_i$ is of the form \eqref{eq2.2} whose elements on the main diagonal are monic polynomials. Here $a_i$, $v_i$ are sizes of the block rows. 
Set 
\[
\cU_2 \!:=\!\m{I & 0 \\ 0 & \tcU_2}, \
\cP_3 \!:=\! \cU_2 \cP_2 \cV_2 \!=\!
 \m{\hcP_{11} & \vline & \hcP_{12} \\ \hline  0 &\vline &  0 \\ 0 &\vline &  0}, \
\cQ_3 \!:=\! \cU_2 \cQ_2 \cV_2 \!=\! 
 \m{\hcQ_{11} & \vline & \hcQ_{12} \\ \hline  0 &\vline &  \Si_i \\ 0 &\vline &  0},
\]
increase $i$ by 1 and repeat the process from Step 1 by applying unimodular transformations for the pair 
$(\hcP_{11},\hcQ_{11})$ with appropriate embedding to the complete matrix polynomial pair.

\noindent \textbf{End.}

Clearly, this procedure terminates after a finite number of iterations since \break $\sum_{i \geq 1} a_i \! \leq \!n.$ 
Thus, we arrive at the matrix polynomial pair $(\tcP,\tcQ )$ where
\[
\tcP := \left[ \begin{array}{cc|c|c|c|c}
\Si_p &0 & * & *  & \dots   &*  \\ \hline
        &   & 0 & *  & \dots   &*  \\
        &   & 0 & *  & \dots   &*  \\ \hline
        &   &    & 0  & \dots   &*  \\
        &   &    & 0  & \dots   &*  \\  \hline              
        &   &    &     & \ddots &\vdots  \\  \hline         
        &   &    &     &           &0        \\              
        &   &    &     &           &0     
\end{array} \right], \ 
\tcQ = \left[ \begin{array}{cc|c|c|c|c}
\cQ_{11} & \cQ_{12} & *        & *            & \dots   & *  \\ \hline
               &                & \Si_q  & *            & \dots   & *  \\
               &                & 0        & *            & \dots   & *  \\ \hline
               &                &           & \Si_{q-1}& \dots   & * \\
               &                &           & 0            & \dots   & *  \\  \hline              
               &                &           &               &\ddots  & \vdots  \\  \hline         
               &                &           &               &           & \Si_1 \\
               &                &           &               &           & 0
\end{array} \right].
\]
Permuting the $4^{th}$, $6^{th}$, $8^{th},\dots$ block rows to the end, we then have 
\eqref{staircase-form}.
\qed

Notice that, due to the appearance of $\cV$, we have to assume that an initial function $x|_{[-\tau,0]}$ is as smooth as needed.
Then, Theorem \ref{thm2.1} applied to the DDAE \eqref{eq1.1} gives us the following result.

\begin{theorem}\label{thm2.2}
Consider the DDAE \eqref{eq2.3} and assume that an initial function $x|_{[-\tau,0]}$ is sufficiently smooth. 
Then system \eqref{eq2.3} is equivalent (in the sense that there is a 
bijective mapping between the solution spaces via a unimodular matrix polynomial) 
to the following system
\begin{equation}\label{eq2.4}
 \m{\Si_{P}  & 0 & \cP_{13} \\ 0 & 0 & \cP_{23} \\ 0 & 0 & \cP_{33}} \m{y_1(t) \\ y_2(t) \\ y_3(t)} 
 = \m{\cQ_{11}  & \cQ_{12} & \cQ_{13} \\ 0 & 0 & \cQ_{23} \\ 0 & 0 & \cQ_{33}} 
 \m{y_1(t-\tau) \\ y_2(t-\tau) \\ y_3(t-\tau)}  + \m{f_1 \\ f_2 \\f_3},
\end{equation}
where $\Si_{P}$ is as in \eqref{staircase-form}, $\cP_{23}$, $\cP_{33}$, $\cQ_{33}$ 
are block upper triangular with zero diagonal blocks, and $\cQ_{23}$ is block upper 
triangular with diagonal blocks of the form \eqref{eq2.2}.
\end{theorem}
\proof
Consider the matrix polynomial pair $(\cP,\cQ)$ associated with the DDAE \eqref{eq2.3}, 
we apply Theorem \ref{thm2.1} to get two unimodular matrix polynomials $\cU$, $\cV$ 
such that $(\cU\cP\cV,\cU\cQ\cV)$ takes the form \eqref{staircase-form}.
Changing the variable $x= \cV(\frac{d}{dt}) y$, and scaling the 
system \eqref{eq2.3} with $\cU$, it immediately leads to \eqref{eq2.4}.
\qed

To illustrate the applicability of Theorem \ref{thm2.1} in the solvability analysis of the DDAE \eqref{eq2.3}, we consider two following examples.

\begin{example}\label{exam2.1}
Consider the following DDAE on the time interval $[0,\infty)$
	\[
	\m{0 & 0 & 1 \\ 0 & 0 & 0 \\ 0 & 0 & 0} \m{\dot{x}_1(t) \\ \dot{x}_2(t) \\ \dot{x}_3(t)} = 
    \m{0 & 1 & 0 \\ 0 & 0 & 1 \\ 0 & -1 & 0} \m{x_1(t) \\ x_2(t) \\ x_3(t)} + 
	\m{0 & 0 & 0\\1 & 0 & 0\\0 & 0 & 0} \m{x_1(t-\tau) \\ x_2(t-\tau) \\ x_3(t-\tau)} + \m{f_1(t) \\ f_2(t) \\ f_3(t)}.
	\]
First we rewrite the system in the matrix polynomial form 
\begin{equation}\label{eq2.6}
 \m{0 & -1 & \frac{d}{dt} \\ 0 & 0 & -1 \\ 0 & 1 & 0} \m{x_1(t) \\ x_2(t) \\ x_3(t)} = \m{0 & 0 & 0\\1 & 0 & 0\\0 & 0 & 0} \m{x_1(t-\tau) \\ x_2(t-\tau) \\ x_3(t-\tau)} + \m{f_1(t) \\ f_2(t) \\ f_3(t)}.
\end{equation}
Theorem \ref{thm2.1} applied to \eqref{eq2.6} terminates after two steps as follows:\\
Step 1: $\cU_1 = \m{0 & 0 & 1 \\ 0 & -1 & 0 \\ 1 & \lb & 1}$, $\cV_1=\m{0 & 0 & 1 \\ 1 & 0 & 0 \\ 0 & 1 & 0}$, 
$\cP_2=\m{1 & 0 & 0 \\ 0 & 1 & 0 \\ 0 & 0 & 0}$, $\cQ_1=\m{0 & 0 & 0 \\ 0 & 0 & -1 \\ 0 & 0 & \lb}$.\\
Step 2: $\cU_2=\cV_2=I_3$, $\cP_3=\cP_2$, $\cQ_3=\cQ_2$.\\
The transformed system \eqref{eq2.4}, therefore, is
\bens
\m{1 & 0 \\ 0 & 1} \m{y_1(t) \\ y_2(t)} &=& \m{0 \\ -1}y_3(t-\tau) + \m{f_3(t) \\ -f_2(t)}, \\
                                     0  &=& \frac{d}{dt} y_3(t-\tau) + f_1(t) + \frac{d}{dt}f_2(t) + f_3(t),
\eens
where $x = \cV_1 \cV_2 y = \m{y_3 & y_1 & y_2}$. Clearly, this system gives us an explicit solution to $y$, and hence, we obtain the solution $x$ to \eqref{eq2.6}.
\end{example}

\begin{example}\label{exam2.2}
Consider the following DDAE on the time interval $[0,\infty)$
	\begin{equation}\label{eq2.7}
	\m{0 & 0 & 0 \\ 0 & 1 & 0 \\ 0 & 0 & 0} \m{\dot{x}_1(t) \\ \dot{x}_2(t) \\ \dot{x}_3(t)} = 
	\m{0 & 1 & 0 \\ 0 & 0 & 1 \\ 0 & 0 & 0} \m{x_1(t) \\ x_2(t) \\ x_3(t)} + 
	\m{0 & -1 & 0 \\ 0 & 0 & 0 \\ 0 & 0 & 1} \m{x_1(t-\tau) \\ x_2(t-\tau) \\ x_3(t-\tau)} + \m{f_1(t) \\ f_2(t) \\ f_3(t)}.
	\end{equation}
Theorem \ref{thm2.1} applied to \eqref{eq2.6} terminates after two steps as follows:\\
Step 1: $\cU_1 = \m{-1 & 0 & 0 \\ -\lb & -1 & 0 \\ 0 & 0 & 1}$, $\cV_1=\m{0 & 0 & 1 \\ 1 & 0 & 0 \\ 0 & 1 & 0}$, 
$\cP_2=\m{1 & 0 & 0 \\ 0 & 1 & 0 \\ 0 & 0 & 0}$, $\cQ_2=\m{1 & 0 & 0 \\ \lb & 0 & 0 \\ 0 & 1 & 0}$.\\
Step 2: $\cU_2=I_3$, $\cV_2=\m{0 & 1 & 0 \\ 0 & 0 & 1 \\ 1 & 0 & 0}$, 
$\cP_3=\m{0 & 1 & 0 \\ 0 & 0 & 1 \\ 0 & 0 & 0}$, $\cQ_3=\m{0 & 1 & 0 \\ 0 & \lb & 0 \\ 0 & 0 & 1}$.\\
We notice that $x = \cV_1 \cV_2 y = y$, and hence, the transformed system \eqref{eq2.4} is
\[
\m{0 & 1 & 0 \\ 0 & 0 & 1 \\ 0 & 0 & 0} \m{x_1(t) \\ x_2(t) \\ x_3(t)} = 
\m{0 & 1 & 0 \\ 0 & \frac{d}{dt} & 0 \\ 0 & 0 & 1} \m{x_1(t-\tau) \\ x_2(t-\tau) \\ x_3(t-\tau)} + \m{-f_1(t) \\ -\dot{f}_1(t)-f_2(t) \\ f_3(t)}.
\]
Clearly, the last equation of this system gives us an explicit solution to $x_3(t)$. Inserting it into the second equation one can solve $x_2$. Finally, the component $x_1$ can be freely chosen, which can be reinterpreted as an input, while the first equation is in fact only the consistency condition for an 
inhomogeneity.
\end{example}

Making use of the transformed system \eqref{eq2.4}, we deduce the solvability properties of the DDAE \eqref{eq2.3} as follows.

\begin{corollary}\label{coro2.0}
Consider the DDAE \eqref{eq2.3} and the transformed system \eqref{eq2.4}. \linebreak Then the 
following claims hold.
\begin{enumerate}
\item[i)] The component $y_3$ is fixed by the equation 
\begin{equation}\label{eq2.5}
 \cP_{23}(\fr{d}{dt}) y_3(t) = \cQ_{23}(\fr{d}{dt}) y_3(t-\tau) + f_2(t). 
\end{equation}
\item[ii)] The third equation of \eqref{eq2.4} provides only consistency conditions for $f_3$.
\item[iii)] The second component $y_2$ can be freely chosen, so it will be reinterpreted as an input.
\end{enumerate}
\end{corollary}
\proof
Rewrite equation \eqref{eq2.5} in the explicit form, we have
\[
 \m{       0  & *   & \dots  & *  \\
              & 0   & \dots  & *  \\
              &     & \ddots & \vdots  \\  
              &     &        & 0       
 } y_3(t) = 
 \m{ \Si_{q}  & *   & \dots  & *  \\
              &\Si_{q-1}& \dots  & *  \\
              &     & \ddots & \vdots  \\  
              &     &        & \Si_{1}
   } y_3(t-\tau) + f_2(t).
\]
Partition $y_3(t)$ conformably as $y_3(t)=\m{z^T_q(t) & z^T_{q-1}(t) & \dots & z^T_1(t)}^T$,
we can recursively solve for the components $z_1$, $z_2,\dots, z_q$, and hence 
$y_3$ is fixed. The remaining assertions are straightforward.
\qed

Another important contribution of the condensed form \eqref{staircase-form} is to point out the relation between the unique solvability of system \eqref{eq2.3} and its spectral property, as stated in the following corollary. 

\begin{corollary}\label{coro2.1} Assume that an initial function associated with the DDAE \eqref{eq2.3} is consistent and sufficiently smooth. 
Then system \eqref{eq2.3} has a unique piecewise solution if and only if the characteristic polynomial $\det(\cP(\lb)- \om \cQ(\lb))$ is not identically zero.
\end{corollary}
\proof
Due to Corollary \ref{coro2.0}, system \eqref{eq2.3} has a unique piecewise solution 
if and only if in system \eqref{eq2.4}, the second block column and the 
last block row do not appear. In this case, the DDAE \eqref{eq2.3} is of square size and 
\[
 \det(\cU) \det(\cP(\lb)- \om \cQ(\lb)) \det(\cV) = 
 \det\left( 
 \m{\Si_p-\om \cQ_{11} &  P_{13}-\om Q_{13} \\ 0 & P_{23} - \om Q_{23}} \right) \ .
\]
Therefore,
\begin{eqnarray*}
 & & \det(\cU) \det(\cP(\lb)- \om \cQ(\lb)) \det(\cV) \\
 &=& \det(\Si_P - \om \cQ_{11}) \cdot 
 \det \left( \m{ - \om \Si_{q} & *              & \dots  & *  \\
                               &- \om \Si_{q-1}   & \dots  & *  \\
                               &                & \ddots & \vdots  \\  
                               &                &        & - \om \Si_{1}       
 }\right) \\
 &=& \det(\Si_P - \om \cQ_{11}) \ \overset{q}{\underset{i=1}{\prod}}\det(-\om \Si_i).
\end{eqnarray*}
Since $\cU$, $\cV$ are unimodular, their determinants are nonzero constants, and hence, we obtain the desired result.
\qed

Applying Corollary \ref{coro2.1} to system \eqref{eq1.1}, the following corollary is evident. 

\begin{corollary}\label{coro2.2} Assume that an initial function $x|_{[-\tau,0]}$ 
is consistent and sufficiently smooth. Then the DDAE \eqref{eq1.1} is uniquely solvable if and only if the matrix triple $(E,A,B)$ is regular.
\end{corollary}

\begin{remark} We notice that the behaviour method presented in this section is useful for analyzing the theoretical solvability of \eqref{eq1.1}. In fact, due to the author's knowledge, till now this is the first time that the solvability of \eqref{eq1.1} can be throughly analyzed in the case that the matrix pair $(E,A)$ is non-regular. Nevertheless, it is not good enough for the numerical computation of the solution, due to the fact that applying Smith form may reduce an index, and therefore, the complexity of the problem itself. Another reason is that, in order to use the behaviour approach, one needs to assume that the unknown function $x$ to be infinitely differentiable, \cite{GohLR82}. 
This obstacle can be overcome under the consideration of the distributional solution concept, but it may happen that the solution to a new system is no longer the classical solution to the original problem, \cite{Tre09}. 
\end{remark}

%%%%%%%%%%%%%%section 3%%%%%%%%%
%\newpage 
\section{Stability analysis of DDAEs}\label{sec 3}

%%%%%%%%%%%%%%%%%%%%%%%%%%%%%
In this section, we study the stability analysis of \eqref{homDDAE}. As usual, we assume that the considered system is regular, i. e., for any consistent initial function $\vphi$, there exists a unique solution $x(t)$. Recall that Corollary \ref{coro2.2} implies that the necessary condition for the regularity of system \eqref{homDDAE} is the non-emptiness of the resolvent set $\rho(E,A,B)$. Now let us recall one important result for linear homogeneous DDEs, taken from \cite{HalL93}.
\begin{proposition}\label{Prop3.0} Consider a linear homogeneous DDE of the form	
\[ \dot{x}(t) = A x(t) + Bx(t-\tau), \ \mbox{ for all } t\in [0,\infty).
 \]	
Then it is exponentially stable if and only if $\si(I,A,B)\subseteq \C_{-}$.
\end{proposition}

In comparison with DDEs, to introduce a new concept of exponential stability for the DDAE \eqref{homDDAE}, the first and most natural idea would be 
adding a consistency assumption on an initial function $\vphi$, see e.g. \cite{Mic11}. We rephrase it in the next definition.
\begin{definition}\label{def3.2}
The null solution $x=0$ of the DDAE \eqref{homDDAE} is called \emph{exponentially stable} if there exist positive constants 
$\de$ and $\ga$ such that for any consistent initial function $\vphi \in C([-\tau,0],\R^n)$, the solution $x=x(t,\vphi)$ of the corresponding 
IVP to \eqref{homDDAE} satisfies
 \[
   \|x(t)\| \leq \de e^{-\ga t} \|\vphi\|_{\infty}, \ \mbox{ for every } \ t\geq 0.
 \]
\end{definition}
For the exponential stability of DDAEs, let us recall one important result presented in \cite{DuLMT13}.
\begin{proposition}\label{Prop3.1} Assume that the DDAE \eqref{homDDAE} has the same solution set as the so-called \emph{strangeness-free formulation}, which takes the form
\be\label{s-free form}
\m{E_1 \\ 0} \dot{x}(t) = \m{A_1 \\ A_2} x(t) + \m{B_1 \\ B_2} x(t-\tau), \ \mbox{ for all } t\in [0,\infty),
\ee
where $\m{E_1 \\ A_2} \in \R^{n,n}$ is nonsingular. Then, \eqref{homDDAE} is exponentially stable if and only if $\si(E,A,B)\subseteq \C_{-}$.
\end{proposition}

However, inherited from DAE theory, the solution $x(t)$ usually depends not only on $x(t-\tau)$ but also on its derivatives $\dot{x}(t-\tau),\dots,x^{(\mu)}(t-\tau)$, for some $\mu\in \N$, which is called the \emph{strangeness-index} of system \eqref{eq1.1}. 
Therefore, Proposition \ref{Prop3.0} is no longer valid for general high-index DDAEs. We demonstrate this fact in the following example.
\begin{example}\label{exam3.1}
Consider the following DDAE on the time interval $[0,\infty)$
 \begin{equation}\label{eq3.5}
  \m{1 & 0 & 0\\0 & 0 & 1\\0 & 0 & 0} \m{\dot{x}_1(t) \\ \dot{x}_2(t) \\ \dot{x}_3(t)} = 
  \m{-1 & 0 & 0\\0 & 1 & 0\\0 & 0 & 1} \m{x_1(t) \\ x_2(t) \\ x_3(t)} + 
  \m{0 & 0 & 0\\0 & 0 & 0\\-1 & 0 & 0} \m{x_1(t-\tau) \\ x_2(t-\tau) \\ x_3(t-\tau)}.
 \end{equation}
Taking derivative of $x_3(t)$ from the third equation, and substituting it into the second one, we obtain a new system
 \begin{subequations}\label{eq3.6}
 \begin{eqnarray}
  \label{eq3.6a} \dot{x}_1(t) &=&-x_1(t), \\
  \label{eq3.6b} 0 &=& x_2(t) - \dot{x}_1(t-\tau), \\
  \label{eq3.6c} 0 &=& x_3(t) - x_1(t-\tau).
 \end{eqnarray} 
 \end{subequations}
Clearly, \eqref{eq3.6b} implies that system \eqref{eq3.5} is not stable, since on the interval $[0,\tau]$ we have $x_2(t)=\dot{\vphi}_1(t-\tau)$. Nevertheless, one can directly verify that the spectrum $\si(E,A,B)$ is $\si(E,A,B)=\{-1\} \subseteq \C_{-}$. 

Besides that, the existence of a solution $x$ is obtained when an initial function $\vphi$ belongs to the space $C^1([-\tau,0],\R^n)$. If this is the case, the solution of system \eqref{eq3.5} is
\begin{subequations}
	\begin{eqnarray*}
		x_1(t)&=&e^{-t}\vphi(0), \ \mbox{ for all } t\in [0,\infty), \\
		x_2(t)&=&
		\begin{cases}
			& \dot{\vphi}_1(t), \ \mbox{ for all } t\in [0,\tau],\\
			&-e^{-(t-\tau)}\vphi(0), \ \mbox{ for all } t\in [\tau,\infty), 			
		\end{cases} \\
		x_3(t)&=&
		\begin{cases}
			& \vphi_1(t-\tau), \ \mbox{ for all } t\in [0,\tau], \\
			& e^{-(t-\tau)}\vphi(0), \ \mbox{ for all } t\in [\tau,\infty).
		\end{cases}
	\end{eqnarray*}  
\end{subequations}
Thus, for the Euclidean norm $\n{\cdot}_2$ and $\vphi \in (C^1([-\tau,0],\R^n),\|\cdot\|_{C^1})$ we obtain the following estimation
\begin{equation*}
\|x(t)\|_2 \leq 3e^{\tau} e^{-t} \|\vphi\|_{C^1}.
\end{equation*}
\end{example}
Example \ref{exam3.1} raises two questions. Firstly, for which type of DDAEs, the condition $\si(E,A,B) \subseteq \C_{-}$ still implies the exponential stability of the system. Secondly, for DDAEs of high-index, how to generalize the stability concept by reducing the phase space of initial functions in such a way that systems like \eqref{eq3.5} are still exponentially stable.

\begin{definition}
The DDAE \eqref{eq1.1} is called \emph{non-advanced} (or \emph{impulse-free}) if for any function $f\in C^{\infty}([0,\infty),\R^n)$ 
and any consistent initial function $\vphi \in C([-\tau,0],\R^n)$, there exists a unique solution $x$ to the IVP \eqref{eq1.1}, \eqref{eq1.4}.
\end{definition}

The following lemma, taken from \cite{HaM16}, gives a strangeness-free formulation for DDAEs.

\begin{lemma}\label{Lem3.1}
Consider the DDAE \eqref{homDDAE}. Furthermore, assume that the IVP \eqref{homDDAE}, \eqref{eq1.4} has a unique solution for every consistent 
initial function $\vphi$. Moreover, assume that the DDAE \eqref{homDDAE} is non-advanced. Then \eqref{homDDAE} can be transformed to the strangeness-free formulation 
\eqref{s-free form}. 
\end{lemma}

Combine Proposition \ref{Prop3.1} and Lemma \ref{Lem3.1}, we obtain the following theorem, which completely characterizes the exponential stability of the DDAE \eqref{homDDAE}.
\begin{theorem}\label{Thm3.2}
Consider the linear, homogeneous DDAE \eqref{homDDAE}. Then, \eqref{homDDAE} is exponentially stable if and only if the following assertions hold.
 \begin{enumerate}
  \item[i)] The DDAE \eqref{homDDAE} is non-advanced. 
  \item[ii)] The spectrum $\si(E,A,B)$ lies entirely on the left half plane.  
 \end{enumerate} 
 \end{theorem}

Now let us move to the second question mentioned above. Example \ref{exam3.1} motivates a new concept of exponential stability for DDAE.
\begin{definition}\label{def3.3}
The null solution $x=0$ of the DDAE \eqref{homDDAE} is called \emph{$C^p$-weakly exponentially stable ($C^p$-w.e.s)} if there exist an integer 
$0\leq p \leq \infty$ and positive constants $\de$ and $\ga$ such that for any consistent initial function $\vphi \in C^p([-\tau,0],\R^n)$, the solution $x=x(t,\vphi)$ of the corresponding IVP for \eqref{homDDAE} satisfies
	\[ \|x(t)\| \leq \de e^{-\ga t} \|\vphi\|_{C^p}, \ \mbox{ for all } t\geq 0. \]
\end{definition}
Clearly, system \eqref{eq3.5} fits perfectly into this definition, where $\de=3e^{\tau}$ and $\gamma=1$. 
Notice that the (classical) exponential stability is exactly $C^0$-w.e.s.. 
Furthermore, even though $C^p$-w.e.s. has been considered for ODEs and PDEs as well, till now we are not aware of any reference for DDAEs.

\begin{theorem}\label{Thm3.3}
Consider the DDAE \eqref{homDDAE} and assume that the matrix triple $(E,A,B)$ is in the block upper triangular form
\be\label{block-triple}
(E,A,B) = \left( \m{I & E_2 \\ 0 & N}, \m{A_1 & A_2 \\ 0 & I}, \m{B_1 & B_2 \\ 0 & B_4} \right), 
\ee
where the matrix $N$ is nilpotent of nilpotentcy index $\nu$. Moreover, suppose that the matrices $N$ and $B_4$ commute. Then the DDAE \eqref{homDDAE} is $C^{\nu}$-w.e.s. if and only if the spectrum $\si(E,A,B)$ 
satisfies $\si(E,A,B) \subseteq \C_{-}$.
\end{theorem}
\proof
Partitioning the variable $x$ appropriately, we can decompose the DDAE \eqref{homDDAE} as follows.
\bens
\dot{x}_1(t) &\!=\!& A_1 x_1(t) \!+\! B_1 x_1(t-\tau) + \left(-E_2\dot{x}_1(t)\!+\!A_2 x_2(t)\!+\!B_2 x_2(t-\tau) \right), \label{eq3.7a} \\
N \dot{x}_2(t) &\!=\!& x_2(t) \!+\! B_4 x_2(t-\tau), \ \mbox{ for all } t\in [0,\infty). \label{eq3.7b}
\eens
Notice that $\si(E,A,B)=\si(I,A_1,B_1) \cup \si(N,I,B_4)$, and hence, the condition that $\si(E,A,B)$ is a subset $\C_{-}$ guarantees that both 
$\si(I,A_1,B_1)$ and $\si(N,I,B_4)$ are also subsets of $\C^{-}$.
As will be seen later in Theorem \ref{Thm4.4}, the commutativity of $N$ and $B_4$, along with the spectral condition $\si(N,I,B_4) \subseteq \C_{-}$ imply that 
the solution $x_2$ to the corresponding IVP of \eqref{eq3.7b} is $C^{\nu}$-w.e.s.. Hence, for any $\vphi \in C^{\nu}([-\tau,0],\R^n)$, the function 
\ $\|-E_2\dot{x}_1(t)\!+\!A_2 x_2(t)\!+\!B_2 x_2(t-\tau)\|$ \	 is bounded from above by an exponentially decreasing function as $t\rar\infty$. Then, the classical result in \cite{BelC63} guarantees that the first component $x_1$ is exponentially stable. This completes the proof.
\qed

In the following example we demonstrate, that without the commutativity of $N$ and $B_4$, the w.e.s of the solution $x$ does not imply the spectral condition $\si(E,A,B) \subseteq \C_{-}$.

\begin{example}
Consider the following DDAE
\be\label{eq3.8}
 \m{0 & 1 \\ 0 & 0} \m{\dot{x}_1(t) \\ \dot{x}_2(t)} =  \m{1 & 0 \\ 0 & 1} \m{{x}_1(t) \\ {x}_2(t)} +  \m{0 & 0 \\ -\ga & 0} \m{{x}_1(t-1) \\ {x}_2(t-1)}, \ \mbox{ for all } t\in [0,\infty).
\ee  
where $\ga \not= 0$ is a real parameter. Direct computation turns out that $N=\m{0 & 1 \\ 0 & 0}$ and $B_4=\m{0 & 0 \\ -\ga & 0}$ do not commute. 
The spectrum of \eqref{eq3.8} is the solution set of the equation $1-\ga \lb e^{-\lb}=0$, or equivalently, $-\lb e^{-\lb}=1/\ga$. Making use of the Lambert $W$ function \cite{Cor96} one can compute $\lb$ from the above equation. It turns out that for $\gamma=0.5$ there are many eigenvalues $\lb$ of \eqref{eq3.8} with positive real parts, see Figure \ref{Fig1}. Thus, $\si(E,A,B) \not\subseteq \C_{-}$. On the other hand, system \eqref{eq3.8} written in details yields that $x_2(t)=\ga x_1(t-\tau)$	and $x_1(t)=\ga \dot{x}_1(t-\tau)$, and hence \eqref{eq3.8} is $C^{\infty}$-w.e.s. for any $|\ga|<1$.
% A very interesting observation here, is that, the main branch W_0 of the Lambert function seems to show the opposite conclusion to the stability, in comparison to DDEs. In details, for DDEs the exponential stability is obtained if and only if the real part of the main branch (and consequently, others branches) is negative, while for our system it is obtained if and only if the main branch (and consequently, others branches) is positive. This is due to the monotonicity of the real parts of the eigenvalues w.r.t the branch numbers.
%
\begin{figure}[h!]
	\centering
	\includegraphics[scale=0.9]{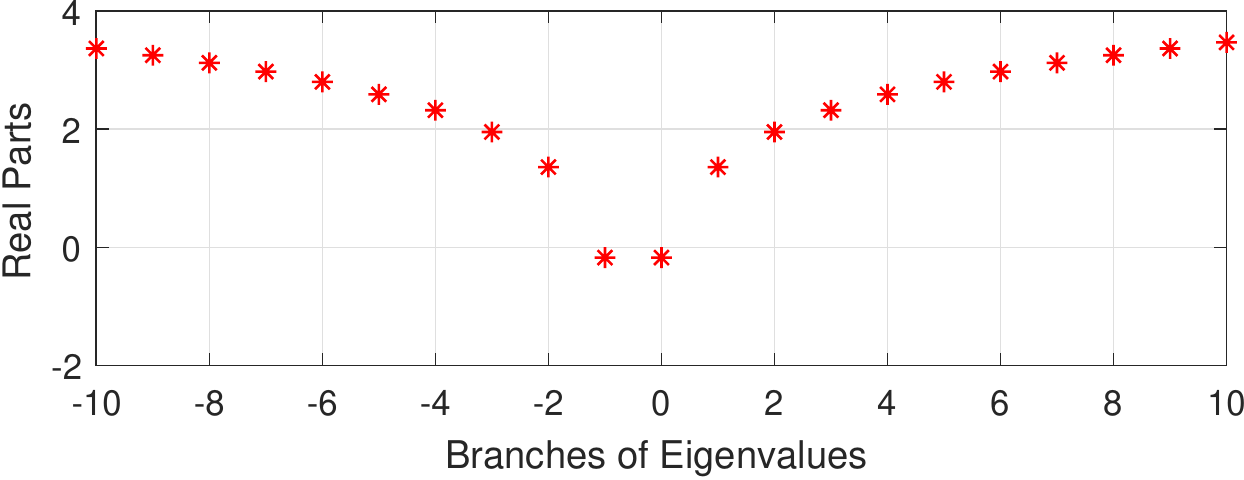}
	\caption{Real part of eigenvalues of \eqref{eq3.8}, where $\gamma=0.5$.}
	\label{Fig1}
\end{figure}
\end{example}

%%%%%%%%%%%%%%section 4%%%%%%%%%
%\newpage 
\section{DDAEs with commutative coefficients}\label{sec4}
This section is devoted to DDAEs whose the matrix coefficients pairwise commute. We will focus on the explicit representation of the solution to the DDAE \eqref{eq1.1} and the stability of the associated homogeneous DDAE \eqref{homDDAE}.
Similar to the case of DAEs, one can transform the DDAE \eqref{eq1.1} by employing \emph{globally equivalent transformations} defined as follows.
\begin{definition} Two triples of matrices $(E_1,A_1,B_1)$ and $(E_2,A_2,B_2)$ in $(\C^{\ell,n})^3$ are called {\it (globally) equivalent } if there exist nonsingular matrices $P\in \C^{\ell,\ell}$ and $Q\in \C^{n,n}$ such that $(E_2,A_2,B_2) = (PE_1Q,PA_1Q,PB_1Q)$. If this is the case, we write $(E_1,A_1,B_1) \sim (E_2,A_2,B_2)$.
\end{definition}
\begin{definition}
	Let $E$, $A$, $B \in \C^{n,n}$. The matrix triple $(E,A,B)$ is called a {\it commutative triple} if $E$, $A$, and $B$ pairwise commute.
\end{definition}
The following lemma will be very useful for our study later. 
\begin{lemma}\label{Lem4} Consider five square matrices $J$, $N$, $\tN$, $S$, $U$ of the same size. Then the following assertions hold.
\begin{enumerate}
\item[i)] If $N$ and $S$ commute and $N$ is nilpotent then $NS$ is nilpotent.
\item[ii)] If $N$ is nilpotent and $U$ is invertible then the matrix $U^{-1}NU$ is nilpotent of the same nilpotency index.
\item[iii)] If both $N$, $\tN$ are nilpotent and they commute then for any scalar $\a$, $\b$ the matrix $\a N+ \b \tN$ is also nilpotent.
\item[iv)] If $J$ is invertible, $N$ is nilpotent and they commute then $J-N$ is invertible and the inverse matrix is given by
\[
(J-N)^{-1} = \overset{\nu-1}{\underset{i=0}{\sum}} (J^{-1})^{i+1} N^i,
\]
where $\nu$ is the nilpotency index of $N$.
\end{enumerate}
\end{lemma}
\proof
The simple proof can be found in \cite{HorJ90}.
\qed

\begin{lemma}\label{Lem4.3}
	Consider a commutative matrix triple $(E,A,B) \in (\C^{n,n})^3$. Then there exists a nonsingular matrix $U\in \C^{n,n}$ such that
	\[
	(UEU^{-1},UAU^{-1},UBU^{-1}) = \left( \m{J & 0 \\ 0 & N},\m{A_{11} & 0 \\ 0 & A_{22}},\m{B_{11} & 0 \\ 0 & B_{22}} \right),
	\]
	where $J$ is nonsingular, $N$ is nilpotent. Furthermore, both triples $(J,A_{11},B_{11})$ and $(N,A_{22},B_{22})$ are commutative.
\end{lemma}
\proof
	First, making use of the Jordan canonical form, we see that there exists a nonsingular matrix $U$ such that
	\[ (UEU^{-1},UAU^{-1},UBU^{-1}) = \left( \m{J & 0 \\ 0 & N},\m{A_{11} & A_{12} \\ A_{21} & A_{22}},\m{B_{11} & B_{12} \\ B_{21} & B_{22}} \right),
	\]
	where $J$ is nonsingular, $N$ is nilpotent. Due to the commutativity of the matrix pair $(E,A)$, the matrices $UEU^{-1}$ and $UAU^{-1}$ also commute, and hence
	\be\label{eq4.17}
	NA_{21} = A_{21} J.
	\ee
	Let $\nu_N$ be the nilpotency index of $N$. Scaling \eqref{eq4.17} with $N^{\nu_N-1}$, one has $0 = N^{\nu_N-1} A_{21}J$, and due to the 
	invertibility of $J$, it follows that \ $0 = N^{\nu_N-1} A_{21}$. 
	Continuing in the same way, eventually one has $0 = NA_{21}$ and hence $A_{21} = 0$. 
	Analogously, one obtains $A_{12}=0$, $B_{21}=0$, $B_{12}=0$ and then the proof is finished.
\qed

In the following theorem, we transform commutative matrix triples into their block diagonal form. 
\begin{theorem}\label{Thm4.1}
	Suppose that $(E,A,B) \in (\C^{n,n})^3$ is a commutative triple. Then, there exists a nonsingular matrix $U\in \C^{n,n}$ such that, 
	\begin{align}\label{eq4.4}
	& (UEU^{-1},UAU^{-1},UBU^{-1})  \notag \\
	&=\left(\ES,\AS,\BS \right),
	\end{align}
	where $J^E$, $J^A$, $J^B$ are nonsingular, $N^E_{2}$, $N^E_{3}$, $N^E_{4}$, $N^A_{3}$, $N^A_{4}$, $N^B_{4}$ are nilpotent. Moreover, if the matrix triple $(E,A,B)$ is regular then the last block row and the last block column are not present.
\end{theorem}
\proof Only within this proof, for notational convenience, we use the superscripts $1$, $2$, $3$ for some matrices. 
Applying Lemma \ref{Lem4.3} to the triple $(E,A,B)$, we find a nonsingular matrix $U_1$ such that
	\[
	(U_1 EU_1^{-1}, U_1AU_1^{-1}, U_1BU_1^{-1}) = \left( \m{J^E & 0 \\ 0 & E^{1}_{2}}, \m{A^{1}_{1} & 0 \\ 0 & A^{1}_{2}}, 
	\m{B^{1}_{1} & 0 \\ 0 & B^{1}_{2}} \right),
	\]
	where $J^E$ is nonsingular, $E^{1}_{2}$ is nilpotent and the triple $\left(E^{1}_{2}, A^{1}_{2}, B^{1}_{2} \right)$ is commutative.\\ 
	Now applying Lemma \ref{Lem4.3} to the triple $\left(A^{1}_{2}, E^{1}_{2}, B^{1}_{2} \right)$, we obtain $U_2$ such that 
	\[
	(U_2 A^{1}_{2} U_2^{-1}, U_2 E^{1}_{2} U_2^{-1}, U_2 B^{1}_{2} U_2^{-1}) 
	= \left( \m{J^A & 0 \\ 0 & A^{2}_{2}}, \m{E^{2}_{1} & 0 \\ 0 & E^{2}_{2}}, \m{B^{2}_{1} & 0 \\ 0 & B^{2}_{2}} \right),   
	\]
	where $J^A$ is nonsingular, $A^{2}_{2}$ is nilpotent and the triple $\left(E^{2}_{2}, A^{2}_{2}, B^{2}_{2} \right)$ is commutative. 
	Furthermore, since $E^{1}_{2}$ is nilpotent, due to Lemma \ref{Lem4} ii) we have that $E^{2}_{1}$ and $E^{2}_{2}$ are also nilpotent.\\
	Finally, applying Lemma \ref{Lem4.3} to the triple $\left(B^{2}_{2}, E^{2}_{2}, A^{2}_{2} \right)$, we obtain $U_3$ such that 
	\[
	(U_3 B^{2}_{2} U_3^{-1}, U_3 E^{2}_{2} U_3^{-1}, U_3 A^{2}_{2} U_3^{-1}) 
	= \left( \m{J^B & 0 \\ 0 & B^{3}_{2}}, \m{E^{3}_{1} & 0 \\ 0 & E^{3}_{2}}, \m{A^{3}_{1} & 0 \\ 0 & A^{3}_{2}} \right),   
	\]
	where $J^B$ is nonsingular, $B^{3}_{2}$ is nilpotent. 
	Since $A^{2}_{2}$, $E^{2}_{2}$ are nilpotent, due to Lemma \ref{Lem4} ii) it follows that $E^{3}_{1}$, $E^{3}_{2}$, $A^{3}_{1}$, $A^{3}_{2}$ are also nilpotent.\\	
	Set $\tU_2 := \m{I & 0 \\ 0 & U_2} \in \C^{n,n}$, $\tU_3 := \m{I & 0 \\ 0 & U_3} \in \C^{n,n}$ and 
	$U := \tU_3 \tU_2 U_1$, we have
	\bens
	&\quad& (UEU^{-1},UAU^{-1},UBU^{-1}) \\
	&=& \left(
	\m{J^E & 0 & 0 & 0 \\ 0 & E^{2}_{1}  & 0 & 0 \\ 0 & 0 & E^{3}_{1} & 0 \\ 0 & 0 & 0 & E^{3}_{2}},
	\m{A^{1}_{1} & 0 & 0 & 0 \\ 0 & J^A  & 0 & 0 \\ 0 & 0 & A^{3}_{1} & 0 \\ 0 & 0 & 0 & A^{3}_{2}}, 
	\m{B^{1}_{1} & 0 & 0 & 0 \\ 0 & B^{2}_{1} & 0 & 0 \\ 0 & 0 & J^B & 0 \\ 0 & 0 & 0 & B^{3}_{2}}
	\right),
	\eens
	%
%	\bens
%	&& (UEU^{-1},UAU^{-1},UBU^{-1}) \\
%	&=& \left(
%	\m{J^E & 0 & 0 & 0 \\ 0 & E^{2}_{1}  & 0 & 0 \\ 0 & 0 & E^{3}_{1} & 0 \\ 0 & 0 & 0 & E^{3}_{2}},
%	\m{A^{1}_{1} & 0 & 0 & 0 \\ 0 & J^A  & 0 & 0 \\ 0 & 0 & A^{3}_{1} & 0 \\ 0 & 0 & 0 & A^{3}_{2}}, 
%	\m{B^{1}_{2} & 0 & 0 & 0 \\ 0 & B^{2}_{1} & 0 & 0 \\ 0 & 0 & J^B & 0 \\ 0 & 0 & 0 & B^{3}_{2}}
%	\right),
%	\eens
%	%
	where $J^E$, $J^A$, $J^B$ are nonsingular and $E^{2}_{1}$, $E^{3}_{1}$, $E^{3}_{2}$, $A^{3}_{1}$, $A^{3}_{2}$, $B^{3}_{2}$ are nilpotent. 
	This is exactly the desired form \eqref{eq4.4}.\\[.2cm] 
	To prove the second claim, we assume that the triple $(E,A,B)$ is regular. Then the triple $(UEU^{-1},UAU^{-1},UBU^{-1})$ is also regular. By direct calculation of the characteristic polynomial, we have
	\bens
	&& \det(\lb E-A - e^{-\lb \tau} B) \\
	&&= \det(\lb J^E - A_{1} - e^{-\lb \tau} B_{1}) \cdot \det(\lb N_{2}^E-J^A - e^{-\lb \tau} B_{2}) \cdot \\
	&\quad& \cdot \det(\lb N_{3}^E-N_{3}^A - e^{-\lb \tau} J^B) \cdot \det(\lb N_{4}^E-N_{4}^A - e^{-\lb \tau} N_{4}^B).
	\eens
	Applying Lemma \ref{Lem4} iii) twice, we see that both matrices $N_{4}^A \!+\! e^{-\lb \tau} N_{4}^B$ and $\lb N_{4}^E \!-\! N_{4}^A \!-\! e^{-\lb \tau} N_{4}^B$ 
	are nilpotent. Thus, $\det(\lb N_{4}^E \!-\! N_{4}^A \!-\! e^{-\lb \tau} N_{4}^B) \!=\! 0$, and hence, $\det(\lb E-A - e^{-\lb \tau} B)=0$. 
	Thus, the presence of the last block row and the last block column destroys the regularity of the triple $(E,A,B)$ and hence, the second claim is proven.
\qed

The decomposition \eqref{eq4.4} is very helpful for analyzing the solvability and stability analysis of DDAEs, as will be studied below. Without loss of generality, let us assume that the matrix triple $(E,A,B)$ is already in the form \eqref{eq4.4}.
Partitioning the variable $x$ and the inhomogeneity $f$ correspondingly, we obtain the following system 
\bsq\label{eq4.5}
\begin{alignat}{2}
J^E \dot{x}_1(t) &= A_{1}x_1(t) + B_{1}x_1(t-\tau) + f_1(t), \label{eq4.5a} \\
N^E_{2} \dot{x}_2(t)  &= J^A x_2(t) + B_{2}x_2(t-\tau) + f_2(t), \label{eq4.5b}\\
N^E_{3} \dot{x}_3(t) &= N^A_{3} x_3(t) + J^B x_3(t-\tau) + f_3(t), \label{eq4.5c} \\
N^E_{4} \dot{x}_4(t) &= N^A_{4} x_4(t) + N^B_{4} x_4(t-\tau) + f_4(t). \label{eq4.5d}
\end{alignat}
\esq
The following corollary gives us a necessary condition for the unique solvability of system \eqref{eq1.1}.

\begin{corollary}\label{Cor3.0}
Consider the corresponding IVP to the DDAE \eqref{eq1.1} and assume that this IVP is uniquely solvable. Then in the decomposition \eqref{eq4.4} the last block row and the last block column must not appear. Consequently, in the system \eqref{eq4.5} the last equation must not present.   
\end{corollary}
\proof
These claims follow directly from Corollary \ref{coro2.2} and Theorem \ref{Thm4.1}.
\qed

Rather than \eqref{eq4.5d}, whose presence destroys the regularity of the triple $(E,A,B)$, now we consider the other equations, and we will give explicit representations for their solutions. \\

Let us denote by $\ddt$ the differentiation operator, which maps a continuously differentiable function $w(t)$ to its first order derivative $\dot{w}(t)$, and, 
by $\De_{-\tau}$ the shift forward operator, which maps a function $w(t)$ to $w(t+\tau)$. 
Scaling the equations \eqref{eq4.5a}-\eqref{eq4.5c} with suitable matrices, we obtain
\bens
\dot{x}_1(t) &=& (J^E)^{-1} A_{1}x_1(t) + (J^E)^{-1} B_{1}x_1(t-\tau) + (J^E)^{-1} f_1(t), \\
(J^A)^{-1} N^E_{2} \dot{x}_2(t)  &=& x_2(t) + (J^A)^{-1}B_{2}x_2(t-\tau) + (J^A)^{-1} f_2(t), \\
(J^B)^{-1}N^E_{3} \dot{x}_3(t) &=& (J^B)^{-1}N^A_{3} x_3(t) + x_3(t-\tau) + (J^B)^{-1} f_3(t),
\eens
which will be rewritten as follows 
\bsq\label{eq4.6}
\begin{alignat}{2}
\dot{x}_1(t)            &= \tA_{1}x_1(t) + \tB_{1}x_1(t-\tau) + \tf_1(t), \ \mbox{ $n_1$ equations} \label{eq4.6a} \\
\tN^E_{2} \dot{x}_2(t)  &= x_2(t) + \tB_{2}x_2(t-\tau) + \tf_2(t),  \ \quad \mbox{ $n_2$ equations} \label{eq4.6b}\\
\tN^E_{3} \dot{x}_3(t)  &= \tN^A_{3} x_3(t) + x_3(t-\tau) + \tf_3(t). \ \quad \mbox{ $n_3$ equations} \label{eq4.6c}
\end{alignat}
\esq
% 
%where $\tA_{1}\!:=\!(J^E)^{-1} A^{1}_{1}$, $\tB_{1}\!:=\!(J^E)^{-1} B^{1}_{2}$, $\tN^E_{2}\!:=\!(J^A)^{-1} E^{2}_{1}$, $\tB^{2}_{1}\!:=\!(J^A)^{-1}B^{2}_{1}$.
Here an initial function $\vphi$ is partitioned correspondingly as $\vphi=\m{\vphi^T_1 & \vphi^T_2 & \vphi^T_3}^T$. 
We notice that due to the commutativity of the matrix triples $(J^E,A_1,B_1)$, $(N^E_2,J^A,B_2)$, $(N^E_3,N^A_3,J^B)$, all of three matrix pairs $(\tA_{1},\tB_{1})$, $(\tN^E_{2},\tB_{2})$ and $(\tN^E_{3},\tN^A_{3})$ must commute.

In the following lemma we give an explicit solution to equation \eqref{eq4.6c}.
\begin{lemma}\label{Lem4.4}
	Consider equation \eqref{eq4.6c}. Furthermore, suppose that $\tf_3$ is sufficiently smooth. Then \eqref{eq4.6c} has the unique solution
	\begin{equation}\label{eq4.7}
	x_3(t) = - \overset{n_3-1}{\underset{i=0}{\sum}} \left( \tN^E_{3} \ddt  - \tN^A_{3} \right)^i \tf_3(t+(i+1)\tau),
	\end{equation}
	for all $t\in [-\tau,\infty)$.
\end{lemma}

\proof First we rewrite \eqref{eq4.6c} in the operator form 
	\[
	x_3(t-\tau) - \tN^E_{3}  \ddt  \De_{-\tau}  x_3(t-\tau) + \tN^A_{3} \De_{-\tau}  x_3(t-\tau) = - \tf_3(t),
	\]
	or equivalently, $\left(I-  {\cal L} \De_{-\tau} \right) x_3(t-\tau) = - \tf_3(t)$, where ${\cal L}:= \tN^E_{3} \ddt  - \tN^A_{3}$. 
	Since the shift forward operator $\De_{-\tau}$ commutes with $\tN^E_{3}$, $\tN^A_{3}$, $\ddt$, then it commutes with ${\cal L}$. 
	Moreover, since $\tN^E_{3}$, $\tN^A_{3}$ are nilpotent matrices of dimension at most $n_3$, and they commute, these imply that ${\cal L}^{n_3}=0$.\\
	This leads us to the solution formula 
	\[
	x_3(t) = - \left(I-  {\cal L} \De_{-\tau} \right)^{-1} \tf_3(t+\tau) = - \overset{n_3-1}{\underset{i=0}{\sum}} {\cal L}^i \De_{-\tau}^i \tf_3(t+\tau),
	\]
	for all $t\geq -\tau$, which is exactly \eqref{eq4.7}. Hence, the proof is completed.
\qed

In the following lemma we give an explicit solution to the equation \eqref{eq4.6b}.
\begin{lemma}\label{Lem4.5}
	Consider the equation \eqref{eq4.6b}. Furthermore, suppose that both the initial function $\vphi_2$ and the inhomogeneity $\tf_2$ are sufficiently smooth. Then \eqref{eq4.6b} has the unique solution
\ben\label{eq4.8}
 	x_2(t) &=& \Big( I - \tN^E_{2} \ddt \Big)^{-1} \big( -\tB_{2} x_2(t-\tau)-\tf_2(t) \big) \notag \\
 	       &=& - \overset{\ze-1}{\underset{i=0}{\sum}} \left( \tN^E_{2} \right) ^{i} \left( \tB_{2} x^{(i)}_2(t-\tau) + \tf^{(i)}_2(t) \right),
\een
for all $t\in [0,\infty)$. Here $\ze$ is the nilpotency index of $\tN^E_{2}$.
\end{lemma}
\proof We rewrite \eqref{eq4.6b} in the operator form 
	\[
      \left( I-\tN^E_{2} \ddt \right) x_2(t) = - \tB_{2} x_2(t-\tau) - \tf_2(t).
	\]
	Applying Lemma \ref{Lem4} iv) we obtain
	\[
	x_2(t) = - \overset{\ze-1}{\underset{i=0}{\sum}} \left( \tN^E_{2} \ddt \right)^i \left(\tB_{2} x_2(t-\tau) + \tf_2(t) \right),
	\]
   which is exactly \eqref{eq4.8}. Hence, the proof is completed.
\qed

Now we consider the delay differential equation \eqref{eq4.6a}. Due to the commutativity of the matrix pair ($\tA_1$,$\tB_1$), an explicit representation of $x_1(t)$ has been established, see e.g. \cite{KhuS03,Pop12}. For the reader's convenience, we recall it here.

\begin{lemma}\label{Lem4.6} Consider the corresponding IVP to the DDE \eqref{eq4.6a} with commutative matrix coefficients $\tA_1$, $\tB_1$. Moreover, assume that an initial function $\vphi_1=x_1|_{[-\tau,0]}$ is continuous. Then, the solution $x_1$ to this IVP has the form:
	\bens 
	x_1(t) = e^{\tA_1 t} e_{\tau}^{\hB_1 (t-\tau)} \vphi_1(0) 
	&+& \int_{-\tau}^{0} e^{\tA_1(t-s)}e_{\tau}^{\hB_1(t-2\tau-s)} \hB_1 \vphi_1(s) ds \\ 
	&+& \int_{0}^{t} e^{\tA_1(t-s)} e_{\tau}^{\hB_1(t-\tau-s)} \tf_1(s) ds,
	\eens
for all $t\geq 0$, where $\hB_1=e^{-\tA_1\tau} \tB_1$ and the matrix $e_{\tau}^{\hB_1(t-\tau-s)}$ is defined via 
%
%\be\label{delay matrix exponential}
\[e^{Dt}_{\tau}=
\begin{cases}
	& I_n, \ \mbox{ for all } -\tau \leq t \leq 0,\\
	& I_n + D t + \cfrac{D^2}{2!}(t-\tau)^2+\dots+\cfrac{D^k}{k!}(t-(k-1)\tau)^k, \\
	& \quad \mbox{ for all } (k-1)\tau \leq t \leq k\tau, \ k=1,\ 2, ...
\end{cases}
\]
Here $e^{Dt}_{\tau}$ is usually called the \emph{delay matrix exponential} associated with $D$.
\end{lemma}

Thus, due to Lemmas \ref{Lem4.4}-\ref{Lem4.6},  we have proven the following theorem. 
\begin{theorem}\label{Thm4.3}
	Consider the corresponding IVP to the DDAE \eqref{eq1.1} and suppose that it is uniquely solvable. Moreover, assume that the following conditions are satisfied.
	\begin{enumerate}
		\item[i)] Both of the initial and the inhomogeneity functions are sufficiently smooth.
		\item[ii)] The matrix triple $(E,A,B)$ is commutative and regular. 
		\item[iii)] The DDAE \eqref{eq1.1} is already in the form \eqref{eq4.5}, which can be rewritten as \eqref{eq4.6}.
	\end{enumerate}
  Then the solution to this IVP has the form 
   \bens
   x_1(t) &=& e^{\tA_1 t} e_{\tau}^{\hB_1 (t-\tau)} \vphi_1(0) + \int_{-\tau}^{0} e^{\tA_1(t-s)}e_{\tau}^{\hB_1(t-2\tau-s)} \hB_1 \vphi_1(s) ds \\
   &\quad& + \int_{0}^{t} e^{\tA_1(t-s)} e_{\tau}^{\hB_1(t-\tau-s)} \tf_1(s) ds,\\
   x_2(t) &=& - \overset{\ze-1}{\underset{i=0}{\sum}} \left( \tN^E_{2} \right) ^{i} \left( \tB_{2} x^{(i)}_2(t-\tau) + \tf^{(i)}_2(t) \right),\\
   x_3(t) &=& - \overset{n_3-1}{\underset{i=0}{\sum}} \left( \tN^E_{3} \ddt  - \tN^A_{3} \right)^i \tf_3(t+(i+1)\tau),
   \eens
 for all $t\geq 0$.
\end{theorem}

In the remaining part of this section we study the stability of the linear homogeneous DDAE \eqref{homDDAE} with commutative coefficients. 
As discussed in Section \ref{sec 3}, the classical stability concept is not always suitable for high index DDAEs. 
The new concept of weakly exponential stability proposed in Definition \ref{def3.3} would be more appropriate. 
Our aim is to analyze the relation between this type of stability for the DDAE \eqref{homDDAE} and the location of its spectrum $\si(E,A,B)$.

We assume again that the triple $(E,A,B)$ is already in the block diagonal form \eqref{eq4.4}. Furthermore, due to Corollary \ref{coro2.2}, the uniqueness of the solution to the corresponding IVP to \eqref{eq1.1} guarantees that the last block row and the last block column do not appear in \eqref{eq4.4}. The spectrum of the triple $(E,A,B)$, therefore, is
\[
\si(E,A,B) = \si \left(
\m{J^E & 0 & 0 \\ 0 & N^E_2  & 0 \\ 0 & 0 & N^E_3},
\m{A_1 & 0 & 0 \\ 0 & J^A  & 0 \\ 0 & 0 & N^A_3}, 
\m{B_1 & 0 & 0 \\ 0 & B_2 & 0 \\ 0 & 0 & J^B}
\right).
\]
%&\!=\!& \Big\{ \lb \in \C \ | \  \det \left(\m{\lb J^E \!-\! A_1 \!-\!  e^{-\tau \lb} B_1 & 0 & 0 \\ 0 & \lb N^E_2 \!-\!  J^A \!-\!  e^{-\tau \lb} B_2  & 0 \\ 0 & 0 & \lb N^E_3 \!-\!  N^A_3 \!-\!  e^{-\tau \lb} J^B } \right) \!= \!0
%\Big\}.
%
Notice that due to Lemma \ref{Lem4} iii) and iv), the matrix $\lb N^E_3\!-\! N^A_3$ is nilpotent and the matrix exponential $\lb N^E_3\!-\! N^A_3 \!-\! e^{-\tau \lb} J^B$ is invertible. Consequently, we have the equality
\bens
\si(E,A,B) &=& \si(J^E,A_1,B_1) \cup \si(N^E_2,J^A,B_2) \notag \\
&=& \si(I,(J^E)^{-1} A_1,(J^E)^{-1} B_1) \cup \si( (J^A)^{-1} N^E_2,I, (J^A)^{-1}B_2), \notag  \\
&=:& \si(I,\tA_1,\tB_1) \cup \si( \tN^E_{2},I, \tB_2). 
\eens
Since we are interested in the stability of homogeneous DDAEs, system \eqref{eq4.6} becomes 
\bsq\label{eq4.10}
\begin{alignat}{2}
\dot{x}_1(t)            &= \tA_{1}x_1(t) + \tB_{1}x_1(t-\tau), \label{eq4.10a} \\
\tN^E_{2} \dot{x}_2(t)  &= x_2(t) + \tB_{2}x_2(t-\tau), \label{eq4.10b}\\
\tN^E_{3} \dot{x}_3(t)  &= \tN^A_{3} x_3(t) + x_3(t-\tau). \label{eq4.10c}
\end{alignat}
\esq
Due to Theorem \ref{Thm4.4}, one sees that the solution $x_3$ of \eqref{eq4.10c} is identically $0$. Thus, the stability of the homogeneous DDAE \eqref{homDDAE} only depends on the stability of two systems \eqref{eq4.10a} and \eqref{eq4.10b}. 

\begin{lemma}\label{Lem4.7}
Consider the spectrum $\si( \tN^E_{2},I, \tB_2)$ of equation \eqref{eq4.10b}. \linebreak If $\si( \tN^E_{2},I, \tB_2) \subseteq \C_{-}$, then $\si(\tB_2)$ lies strictly inside the unit circle. Consequently, there exists a norm $\n{\cdot}$ such that $\n{\tB_2} < 1$.
\end{lemma}
\proof
The first claim of this lemma follows directly from the so-called \emph{L-property} of commutative matrices, \cite{MotT52}.
The second claim of this lemma is Lemma B7, \cite{StuH96}.
\qed

\begin{lemma}\label{Lem4.8}
Consider an element $X$ of some ring equipped with two binary operations denoted by $(\cZ,+,\cdot)$. If there exists $\ze\in \N$ such that $X^\ze=0$ then $(1_{\cZ}-X)^{-1}=1_{\cZ}+X+\dots+X^{\ze-1}$. Furthermore, for each $k\in \N$ the element $((1_{\cZ}-X)^{-1})^{k}$ is a polynomial of the variable $k$ whose degree is at most $\ze-1$. The coefficients of this polynomial are of the form $X^j$, where $j\leq \ze-1$. 
\end{lemma}
\proof
The first claim is evident, due to the fact that 
\[
(1_{\cZ}-X) \cdot (1_{\cZ}+X+\dots+X^{\ze-1}) = 1_{\cZ}-X^\ze = 1_{\cZ}.
\] 
Now we prove the second claim by proving that $\left( 1_{\cZ}+X+\dots+X^{\ze-1} \right)^k$ is a polynomial of $k$ with the degree is at most $\ze-1$.\\
Due to the Multinomial Theorem, we see that 
\bens
\left( 1_{\cZ}+X+\dots+X^{\ze-1} \right)^k &=& \underset{i_0+i_1+\dots+i_{\ze-1}=k}{\sum} \binom{k}{i_0,i_1,\dots,i_{\ze-1}} \prod_{j=1}^{\ze-1} (X^j)^{i_j},  \\ 
&=&  \underset{i_0+i_1+\dots+i_{\ze-1}=k}{\sum} \binom{k}{i_0,i_1,\dots,i_{\ze-1}} X^{ \sum_{j=1}^{\ze-1} j i_j}, 
\eens
where $i_j$, $j=0,\dots,\ze-1$ are nonnegative integers and 
\[
\binom{k}{i_0,i_1,\dots,i_{\ze-1}} = \cfrac{k!}{i_0! \ i_1! \ i_2! \dots  \ i_{\ze-1}!} \ .
\]
Notice that for any $i_0 \!\leq\! k\!-\!\ze$ then $\sum_{j=1}^{\ze-1} i_j \!=\! k\!-\!i_0 \!\geq\! \ze$, and hence $\sum_{j=1}^{\ze-1} j i_j \!\geq\! \ze$. 
Since $X^\ze = 0$, in this case we see that 
\[ X^{ \sum_{j=1}^{\ze-1} j i_j} = 0.
\] 
Thus, we have  
\be\label{eq4.13}
\left( 1_{\cZ}+X+\dots+X^{\ze-1} \right)^k =  \underset{ \underset{i_0 \geq k-\ze+1}{i_0+i_1+\dots+i_{\ze-1}=k}}{\sum} \binom{k}{i_0,i_1,\dots,i_{\ze-1}} X^{ \sum_{j=1}^{\ze-1} j i_j}. 
\ee
Clearly, the amount of tuples $(i_1,\dots,i_{\ze-1})$ that satisfies $\sum_{j=1}^{\ze-1} i_j \leq \ze-1$ is finite. 
Moreover, for any fixed tuple $(i_0,i_1,\dots,i_{\ze-1})$ that satisfies $i_0 \geq k-\ze+1$, then 
\[
\binom{k}{i_0,i_1,\dots,i_{\ze-1}} = \cfrac{k!}{i_0! \ i_1! \ i_2! \dots  \ i_{\ze-1}!} = \cfrac{k \ (k-1) \dots (i_0+1)}{i_1! \ i_2! \dots  \ i_{\ze-1}!}
\]
is a polynomial (of variable $k$) whose order is at most $\ze-1$.
Consequently, the right hand side of \eqref{eq4.13} is a polynomial of the variable k, whose degree is at most $\ze-1$.
Finally, we see that the coefficients of this polynomial are of the form $X^m$, which are non-zero if and only if $m\leq \ze-1$. Thus, the proof is completed.
\qed

The following theorem characterizes the stability of the DDAE \eqref{homDDAE}. 
 
\begin{theorem}\label{Thm4.4} Assuming that the corresponding IVP to the DDAE \eqref{eq1.1} is \linebreak uniquely solvable. Moreover, suppose that the matrix triple $(E,A,B)$ is commutative, and an initial function is sufficiently smooth. Then the following assertions hold.
\begin{enumerate}
\item[i)] The solution is exponentially stable if $\si(E,A,B) \subseteq \C_{-}$, and the matrix $N^E_2$ in the block diagonal form \eqref{eq4.4} is identically $0$. 
\item[ii)] The solution is $C^{\ze}$-w.e.s. if $\si(E,A,B) \subseteq \C_{-}$, where $\ze$ is the nilpotency index of $N^E_2$. 
\end{enumerate}
\end{theorem}

\proof First we see that due to the equality \eqref{eq4.11}, both sets $\si(I,\tA_{1},\tB_{1})$ and $\si( \tN_{2},I, \tB_{2})$ belongs to $\C_{-}$. Proposition \ref{Prop3.0} implies the exponential stability of \eqref{eq4.10a}, and hence, we only need to care about \eqref{eq4.10b}. We recall that $\tN^E_2=(J^A)^{-1} N^E_2$, and hence, the  nilpotency indices of $\tN^E_2$ and $N^E_2$ are equal.\\
i) If $N^E_2=0$, then $\tN^E_2=0$ and \eqref{eq4.10b} becomes $0=x_2(t) + \tB_{2}x_2(t-\tau)$. Lemma \ref{Lem4.7} follows that this difference equation is exponentially stable.\\
ii) Now we will prove the second claim. Making use of Lemma \ref{Lem4.5}, we see that
\[
x_2(t) =  \Big( I - \tN^E_{2} \ddt \Big)^{-1} \big( -\tB_{2}  x_2(t-\tau) \big), \ \mbox{for all} \ t\geq 0.
\]
Therefore, simple induction gives us the explicit representation of $x_2$ in terms of $\vphi_2$ as follows
\[
%x_2(t) =  \Bigg[ \left( I + \tN^E_{2} \ddt + \dots + \left( \tN^E_{2} \ddt \right)^{\ze-1} \right) (-\tB_{2}) \Bigg]^{k} {\vphi}_2(t-k\tau).
x_2(t) =  \Bigg[  \Big( I - \tN^E_{2} \ddt \Big)^{-1} (-\tB_{2}) \Bigg]^{k} {\vphi}_2(t-k\tau), \ \mbox{for all} \ t\geq 0,
\]
where $k:= \floor*{\frac{t}{\tau}}+1$. 
Due to the commutativity of $\tN^E_{2}$ and $\tB_{2}$ we then have
\[
x_2(t) =  (-\tB_{2})^{k}  \Bigg[  \Big( I - \tN^E_{2} \ddt \Big)^{-1} \Bigg]^{k} {\vphi}_2(t-k\tau).
\]
Now applying Lemma \ref{Lem4.8} for $X=\tN^E_{2} \ddt$ and notice that $\zeta$ is exactly the nilpotency index of $\tN^E_{2}$, 
we see that $\Bigg[  \Big( I - \tN^E_{2} \ddt \Big)^{-1} \Bigg]^{k}$ is a matrix polynomial of variable $k$ whose degree is at most 
$\ze-1$. Furthermore, also due to Lemma \ref{Lem4.8}, the coefficients of this polynomial are operators of the form $\left(\tN^E_{2} \ddt \right)^m$, where $m\leq \ze-1$. Let us describe this polynomial as follows 
\[
 \frQ(k, \tN^E_{2}\ddt) :=  \sum_{j=0}^{\ze-1} \left( \tN^E_{2} \ddt \right)^{m_j} k^j, 
\]
where $0\leq m_j \leq \ze-1$ for all $j=1,\dots,\ze-1$. \\
Consequently, for any $s\in [-\tau,0]$ and any $k\in \N$ we have that
 \bens
 x_2(s+k\tau) &=& (-\tB_{2})^{k} \sum_{j=0}^{\ze-1} \left( \tN^E_{2} \ddt \right)^{m_j} k^j   {\vphi}_2(s), \\
        &=&(-\tB_{2})^{k} \sum_{j=0}^{\ze-1} \left( \tN^E_{2} \right)^{m_j} k^j  {\vphi}^{(m_j)}_2(s). 
 \eens
Thus, for any given norm $\n{\cdot}$, let $C:= \max\{ \|(\tN^E_{2})^j \|, \ j=0,\dots,\ze-1 \}$, we then have the following estimation
\[ 
\n{x_2(s+k\tau)} \leq C \ \| \tB_2 \|^{k} \  \sum_{j=0}^{\ze-1}k^j \  \n{\vphi_2}_{C^{\ze-1}} \ \mbox{for all} \ t\geq 0. 
\]
Making use of the second claim of Lemma \ref{Lem4.7}, we can choose a suitable norm such that  $\|\tB_2\|<1$, which follows that $\n{x_2(t)}$ converges exponentially to $0$ as $t\rar +\infty$. This completes our proof. 
\qed  

In the following example, we illustrate our result. 
\begin{example}\label{exam4.1}
Consider the following DDAE on the time interval $[0,\infty)$
\begin{equation}\label{eq4.11} 
\m{ 2 & -4  &-8 \\ -8 & -4 &  2 \\ 12 & 16 & 12} \m{\dot{x}_1(t) \\ \dot{x}_2(t) \\ \dot{x}_3(t)} \!=\! 
\m{28 & 36 & 36 \\ -12 & -14  & -24 \\ -12 & -24 & -14} \m{x_1(t) \\ x_2(t) \\ x_3(t)} \!+\! 
\m{2 & -6 & -6 \\ 2 & 9 &  4 \\ 2 & 4 & 9} \m{x_1(t \!-\! \tau) \\ x_2(t \!-\! \tau) \\ x_3(t \!-\! \tau)}.
\end{equation}
As in Theorem \ref{Thm4.1}, we can find the matrix $U = \m{1 & 2 & 2 \\ 2 &  1 & 2 \\ 2 & 2 & 1}$ and hence, we can transform system \eqref{eq4.11} to the block diagonal form \eqref{eq4.4} which reads
\[
\m{ 10 & 0  & 0 \\ 0 & 0 & 10 \\ 0 & 0 & 0} \m{\dot{y}_1(t) \\ \dot{y}_2(t) \\ \dot{y}_3(t)} = 
\m{ -20 & 0 & 0 \\ 0 & 10  & 0 \\ 0 & 0 & 10} \m{y_1(t) \\ y_2(t) \\ y_3(t)} + 
\m{ 10 & 0 & 0 \\ 0 & 5 &  0 \\ 0 & 0 & 5} \m{y_1(t-\tau) \\ y_2(t-\tau) \\ y_3(t-\tau)}.
\]
where the new variable is $y = U^{-1}x$. Due to Theorem \ref{Thm4.4}, \eqref{eq4.11} is $C^{2}$-w.e.s, even though it is unstable in the classical sense.
\end{example}

\begin{remark}\label{rem4.1}
a) Theorem \ref{Thm4.4}b only provides a sufficient conditions for the $C^p$-w.e.s of DDAEs. Nevertheless, it is not a necessary condition, as can be directly seen from Example \ref{exam3.1}, whose the solution is $C^1$-w.e.s, even though the matrix coefficients do not commute. \\
b) There are certain DDAE systems, where the matrix coefficients do not form a commutative triple. However, after applying global equivalent transformation, then the new coefficients are pairwise commutative. About this issue, we refer the readers to \cite{Ha15}.\\
c) In order to determine $\zeta$ (the nilpotency index of $N^E_2$), the task of computing a matrix $U$ is not alway necessary. We notice that due to Lemma \ref{Lem4.3} and Theorem \ref{Thm4.1}, U is computed based on the Jordan canonical form of certain matrices, which is an unstable problem. One way to overcome this is to use global unitary transformations. However, the price to pay is that we do not have a block diagonal triple as in \eqref{eq4.4}, but a block upper triangular triple. For details, also see \cite{Ha15}.
\end{remark}

%%%%%%%%%%%%%%section 5%%%%%%%%%
\section{Conclusion}\label{conclusion}

%%%%%%%%%%%%%%%%%%%%%%%%%%%%%

In this paper, we characterized the solvability and stability analysis for general linear delay differential-algebraic equations (DDAEs) in terms of spectral conditions. We showed that, like DAEs, solvability properties of DDAEs are closely related to the regularity of their matrix triple. 
However, for the stability one needs to be careful, since the classical eigenvalue-based approach is only valid for non-advanced DDAE systems.
Therefore, rather than using only the classical concept of exponential stability, we proposed and considered a more general concept of $C^p$-weakly exponential stability ($C^p$-w.e.s.) for DDAEs. Then, we discussed the solvability and the weakly exponential stability for a class of DDAEs with pairwise commutative coefficients.  

\bigskip
{\bf Acknowledgment.} The author would like to thank an anonymous referee for his constructive comments and suggestions that improve the quality of this paper. The author also thanks Stephan Trenn for helpful comments and fruitful discussions about the first topic of this article.

%%%%%%%%%%%%%%%%%%%%%%%%%%%%%%%%%%%%%%%%%%%%%%%%%%%%%%%%%%%%%
\bibliographystyle{abbrv}
%\bibliography{Phi_Jan_2018}
%\bibliography{Ha18}
%\include{appendix}

\end{document}